\documentclass[a4paper,12pt]{article}
\usepackage{comment}
\usepackage{cite}
\usepackage{amsmath}
\usepackage{amssymb}
\usepackage{amsfonts}
\usepackage[T1]{fontenc}
\usepackage[utf8]{inputenc}
\usepackage{graphicx}
\usepackage{fancyhdr}
\usepackage{float}
\usepackage{xcolor}
\usepackage{authblk}
\usepackage{mathrsfs}
\usepackage{empheq}
\usepackage[hyphens]{url}
\usepackage{hyperref} 
\usepackage[]{breakurl}

\usepackage{geometry}
\begin{document}

\title{A prime sum involving Bernoulli numbers}

\author[$\dagger$]{Jean-Christophe {\sc Pain}$^{1,2,}$\footnote{jean-christophe.pain@cea.fr}\\
\small
$^1$CEA, DAM, DIF, F-91297 Arpajon, France\\
$^2$Universit\'e Paris-Saclay, CEA, Laboratoire Mati\`ere en Conditions Extr\^emes,\\ 
91680 Bruy\`eres-le-Ch\^atel, France
}

\maketitle

\begin{abstract}
In this note, we propose simple identities for primes, which involve two finite nested sums and Bernoulli numbers. The summations can also be expressed in terms of Bernoulli polynomials.
\end{abstract}

\section{Introduction}

Many works are devoted to the determination of explicit formulas for prime numbers (see the non-exhaustive list of references \cite{Willans1954,Veshenevskiy1962,Gandhi1971,Eynden1974,Golomb1976,Ruiz2005,Atanassov2021}). As an example, the Willans formula \cite{Willans1954} for the $n^{th}$ prime $p_n$ reads
\begin{equation}\label{wil}
p_{n}=1+\sum_{i=1}^{2^{n}}\left\lfloor \left\lfloor\frac{n}{\sum_{j=1}^{i}\left\lfloor \cos^{2}\left(\frac{(j-1)!+1}{j}\pi\right)\right\rfloor}\right\rfloor^{1/n}\right\rfloor,
\end{equation}
where $\lfloor x\rfloor$ denotes the integer part of $x$. Equation (\ref{wil}) can be expressed with the prime-counting function $\pi(x)$ (the function counting the number of prime numbers less than or equal to some real number $x$) as
\begin{equation}
p_n=1+\sum_{m=1}^{2^n}\left\lfloor\left\lfloor\frac{n}{1+\pi(m)}\right\rfloor^{1/n}\right\rfloor.
\end{equation}
Willans also gave the following expression for $\pi(m)$ \cite{Ribenboim1971}:
\begin{equation}
\pi(m)=\frac{\sin^2\left(\pi\frac{\left\{(j-1)!\right\}^2}{j}\right)}{\sin^2\left(\frac{\pi}{j}\right)}
\end{equation}
and Min\'a$\mathrm{\check{c}}$ obtained a formula which does not involve any trigonometric function:
\begin{equation}
\pi(m)=\sum_{j=2}^m\left\lfloor\frac{(j-1)!+1}{j}-\left\lfloor\frac{(j-1)!}{j}\right\rfloor\right\rfloor.
\end{equation}
Gandhi published the expression \cite{Gandhi1971,Golomb1974}:
\begin{equation}
p_n=\left\lfloor 1-\frac{1}{\ln 2}\ln\left( -\frac{1}{2}+\sum_{d/\mathscr{P}_{n-1}}\frac{\mu(d)}{2^d-1}\right)\right\rfloor,
\end{equation}
where $\mu(n)$ is the M\"obius function and $\mathscr{P}_n=p_1\cdots p_{n-1}$. Vassilev-Missana \cite{Vassilev-Missana2001}, following Atanassov, \cite{Atanassov2001} found the following three formulas:
\begin{equation}
p_n=\sum_{j=0}^{C(n)}\left\lfloor\frac{1}{1+\left\lfloor\frac{\pi(j)}{n}\right\rfloor}\right\rfloor
\end{equation}
as well as
\begin{equation}
p_n=-2\sum_{j=0}^{C(n)}\zeta\left(-2\left\lfloor\frac{\pi(j)}{n}\right\rfloor\right)
\end{equation}
where $\zeta(n)$ is the Riemann zeta function and
\begin{equation}
p_n=\sum_{j=0}^{C(n)}\frac{1}{\Gamma\left(1-\left\lfloor\frac{\pi(j)}{n}\right\rfloor\right)}
\end{equation}
with
\begin{equation}
C(n)=\left\lfloor\frac{n^2+3n+4}{4}\right\rfloor.
\end{equation}
Dimitrov performed a detailed comparison of the computation time of different existing formulae for the $n^{th}$ prime number \cite{Dimitrov2019}.

Beyond the interest in explicit formulas for the $n^{th}$ prime number, much effort is devoted to the search for identities (and in particular summations \cite{Glaisher1891,Lionnais1983}). Some of them involve the Mangoldt function \cite{Hardy1979,Berndt1994}. One can also mention the rapidly converging series for the Mertens constant $B_1$ (also known as the Hadamard-de la Vallee-Poussin constant or prime reciprocal constant \cite{Mertens1874,Bach1996,Flajolet1996}):
\begin{equation}
B_1=\gamma+\sum_{k=1}^{\infty}\left[\ln\left(1-\frac{1}{p_k}\right)+\frac{1}{p_k}\right]=\gamma+\sum_{m=2}^{\infty}\frac{\mu(m)}{m}\ln\left[\zeta(m)\right],
\end{equation}
where $\gamma$ is the Euler-Mascheroni constant.

In section \ref{sec2}, we derive simple summations for primes using integer parts and binomial expansions. The final results involve a double sum over Bernoulli numbers. The latter double sum can also be expressed in terms of Bernoulli polynomials.

\section{A family of sums for primes involving Bernoulli numbers}\label{sec2}

Let $p$ be a prime, and let us consider the sum
\begin{equation}\label{def}
\mathscr{S}_q(p)=\sum_{k=1}^{p-1}\left\lfloor\frac{k^{2q+1}}{p}\right\rfloor,
\end{equation}
for $q\geq 1$. Making the change of index $p\rightarrow p-1$, one finds
\begin{equation}
\mathscr{S}_q(p)=\sum_{k=1}^{p-1}\left\lfloor\frac{(p-k)^{2q+1}}{p}\right\rfloor.
\end{equation}
Using the binomial expansion:
\begin{equation}
(p-k)^{2q+1}=\sum_{r=0}^{2q+1}\binom{2q+1}{r}(-1)^rk^rp^{2q+1-r}
\end{equation}
and separating the term $r=q$ from the rest of the summation, one gets
\begin{eqnarray}
\mathscr{S}_q(p)&=&\sum_{k=1}^{p-1}\left(\sum_{r=0}^{2q}\left\lfloor(-1)^r\binom{2q+1}{r}k^rp^{2q-r}\right\rfloor+\left\lfloor-\frac{k^{2q+1}}{p}\right\rfloor\right)\nonumber\\
& &\sum_{k=1}^{p-1}\left(\sum_{r=0}^{2q}(-1)^r\binom{2q+1}{r}k^rp^{2q-r}+\left\lfloor-\frac{k^{2q+1}}{p}\right\rfloor\right).
\end{eqnarray}
Using the definition (\ref{def}), we find
\begin{equation}
2\mathscr{S}_q(p)=\sum_{k=1}^{p-1}\left(\sum_{r=0}^{2q}(-1)^r\binom{2q+1}{r}k^rp^{2q-r}+\left\lfloor\frac{k^{2q+1}}{p}\right\rfloor+\left\lfloor-\frac{k^{2q+1}}{p}\right\rfloor\right).
\end{equation}
Since $p$ is prime and greater than $k$, $k^{2q+1}/p$ is not an integer. Therefore
\begin{equation}
\left\lfloor\frac{k^{2q+1}}{p}\right\rfloor+\left\lfloor-\frac{k^{2q+1}}{p}\right\rfloor=-1,
\end{equation}
yielding
\begin{equation}
2\mathscr{S}_q(p)=\sum_{k=1}^{p-1}\sum_{r=0}^{2q}(-1)^r\binom{2q+1}{r}k^rp^{2q-r}-(p-1)
\end{equation}
\emph{i.e.}
\begin{equation}
2\mathscr{S}_q(p)=\sum_{r=0}^{2q}(-1)^r\binom{2q+1}{r}p^{2q-r}\sum_{k=1}^{p-1}k^r-(p-1).
\end{equation}
One has also, for $r>0$:
\begin{equation}
\sum_{k=1}^{p-1}k^r=\frac{1}{r+1}\sum_{l=0}^r\binom{r+1}{l}B_l~p^{r-l+1},
\end{equation}
where $B_l$ is the Bernoulli number of order $l$ \cite{Carlitz1968,Khan1981,Krishnapriyan1995}. For $r=0$, the sum is equal to $(p-1)$. This enables us to write (isolating the $r=0$ term):
\begin{equation}
\mathscr{S}_q(p)=-\frac{(p-1)}{2}+\frac{1}{2}\left(p^{2q}(p-1)+\sum_{r=1}^{2q}(-1)^r\binom{2q+1}{r}p^{2q-r}\frac{(-1)^r}{r+1}\sum_{l=0}^{r}\binom{r+1}{l}B_l~p^{r-l+1}\right)
\end{equation}
and finally
\begin{empheq}[box=\fbox]{align}\label{final}
\mathscr{S}_q(p)=\frac{(p-1)(p^{2q}-1)}{2}+\frac{1}{2}\sum_{r=1}^{2q}\frac{(-1)^r}{r+1}\binom{2q+1}{r}\sum_{l=0}^{r}\binom{r+1}{l}B_l~p^{2q+1-l},
\end{empheq}
which is the main result of the present work. Using the definition of Bernoulli polynomials:
\begin{equation}
B_n(x)=\sum_{k=0}^n\binom{n}{k}B_k~x^{n-k},
\end{equation}
one can show that
\begin{equation}
\mathscr{S}_q(p)=\frac{B_{2q+2}(-p)+B_{2q+2}(p)-2B_{2q+2}}{2p(2q+2)}-\frac{(p^{2q}+p-1)}{2}.
\end{equation}
Using $B_{2q+2}(-p)=B_{2q+2}(p+1)$, one has
\begin{empheq}[box=\fbox]{align}
\mathscr{S}_q(p)=\frac{B_{2q+2}(p+1)+B_{2q+2}(p)-2B_{2q+2}}{2p(2q+2)}-\frac{(p^{2q}+p-1)}{2}
\end{empheq}
and since $B_{2q+2}(1)=B_{2q+2}(0)=B_{2q+2}$, one can write
\begin{equation}
\mathscr{S}_q(p)=\frac{\left(B_{2q+2}(p+1)-B_{2q+2}(1)\right)+\left(B_{2q+2}(p)-B_{2q+2}(0)\right)}{2p(2q+2)}-\frac{(p^{2q}+p-1)}{2}.
\end{equation}
As an example, in the case $q=1$, one recovers the result given in Refs. \cite{Mathe,Doster1993}:

\begin{equation}
\mathscr{S}_1(p)=\sum_{k=1}^{p-1}\left\lfloor\frac{k^3}{p}\right\rfloor=\frac{(p-1)(p-2)(p+1)}{4}
\end{equation}
and in the cases $q=2, 3$ and 4, one obtains respectively:

\begin{equation}
\mathscr{S}_2(p)=\sum_{k=1}^{p-1}\left\lfloor\frac{k^5}{p}\right\rfloor=\frac{1}{12}(p-2)(p-1) (p+1)(2p^2-2 p+3),
\end{equation}
\begin{equation}
\mathscr{S}_3(p)=\sum_{k=1}^{p-1}\left\lfloor\frac{k^7}{p}\right\rfloor=\frac{1}{24}(p-2)(p-1)(p+1)(3p^4-6p^3+5 p^2-2 p+6)
\end{equation}
and
\begin{equation}
\mathscr{S}_4(p)=\sum_{k=1}^{p-1}\left\lfloor\frac{k^9}{p}\right\rfloor=\frac{1}{20}(p-2)(p-1)(p+1)(2p^6-6p^5+5p^4+3p^2-4 p+5).
\end{equation}
\section{Conclusion}

We obtained, using the binomial expansion and properties of the integer part, simple summations for primes, consisting of two finite nested sums and involving Bernoulli numbers. The new identities can also be expressed without any explicit summation, using Bernoulli polynomials. Such identities may be useful for the determination of primality tests.

\section*{Acknowledgements}

I would like to thank Christophe Vignat (Tulane University and Université Paris-Saclay) for kindly bringing to my attention the fact that the double sum in the right-hand side of Eq. (\ref{final}) can be formulated in terms of Bernoulli polynomials.

\end{document}